\documentclass[12pt]{amsart2000}
\usepackage{amsfonts,amscd2000,amsmath2000,enumerate,verbatim,calc}

\usepackage{times} 
\usepackage{mathptm}
\swapnumbers
\theoremstyle{plain}
\newtheorem{Theorem}{Theorem}[section]
\newtheorem*{main}{Theorem}
\newtheorem*{flenner}{Theorem}
\newtheorem{Lemma}[Theorem]{Lemma}
 
\newtheorem{Proposition}[Theorem]{Proposition}

\theoremstyle{definition}
\newtheorem{Definition}[Theorem]{Definition}
\newtheorem{Example}[Theorem]{Example}

\newtheorem{hilb}[Theorem]{Hilbert functions of bigraded modules}
\newtheorem{s2}[Theorem]{The $(S_2)$ property of Serre}
\newtheorem{jmult}[Theorem]{$\mathbf j$-multiplicities}
\newtheorem{genm}[Theorem]{Generalized Samuel multiplicity}
\newtheorem{notation}[Theorem]{Notation}
\newtheorem{coef}[Theorem]{First coefficient ideals}

\newtheorem*{chunk*}{}
\newtheorem*{ack}{\sc Acknowledgement}
\numberwithin{equation}{Theorem}

\theoremstyle{remark}
\newtheorem{Remark}[Theorem]{Remark}

\newcommand{\into}{\operatorname{\hookrightarrow}}

\newcommand{\depth}{\operatorname{depth}}

\newcommand{\het}{\operatorname{ht}}

\newcommand{\length}{\operatorname{\lambda}}

\newcommand{\Ass}{\operatorname{Ass}}
\newcommand{\Hom}[3]{\operatorname{Hom}_{#1}(#2,#3){}}

\newcommand{\Ann}{\operatorname{Ann}}
\newcommand{\Spec}{\operatorname{Spec}}
\newcommand{\fm}{{\mathfrak m}}

\newcommand{\fp}{{\mathfrak p}}

\newcommand{\ringA}{\text{$(A,\fm)$}}

\begin{document}
\title[A numerical characterization]
{A numerical characterization of the $S_2$-ification of a Rees
  algebra}
\author[C\u{a}t\u{a}lin Ciuperc\u{a}]{C\u{a}t\u{a}lin Ciuperc\u{a}}

\address{Department of Mathematics, University of California, Riverside, CA 92521}
\email{ciuperca@math.ucr.edu}
\begin{abstract} 
Let $A$ be a local ring with maximal ideal $\fm$. For an arbitrary
ideal $I$ of $A$, we define the generalized Hilbert coefficients
$j_k(I)\in \mathbb{Z}^{k+1}$ ($0 \leq k \leq \dim A$). When the ideal $I$
is $\fm$-primary, $j_k(I)=(0,\ldots,0,(-1)^ke_k(I))$, where   $e_k(I)$ is
the classical $k^{\mathrm{th}}$ Hilbert coefficient of $I$. 
Using these
  coefficients we give a numerical characterization of the
  homogeneous components of the $S_2$-ification of $
  S=A[It,t^{-1}]$, extending to not necessarily $\fm$-primary ideals 
  the results obtained in \cite{C}.       
\end{abstract}
\maketitle
\bigskip
\section*{Introduction}
Let $\ringA$ be a formally equidimensional local ring and let $I
\subseteq J$ be two ideals of $A$. When $I$ is $\fm$-primary, Rees
proved that $J$ is contained in the integral closure $\overline I$ of
$I$ if and only if $I$ and $J$ have
the same multiplicity. B\"{o}ger \cite{B} extended this result as follows: \emph{let $I
\subseteq J \subseteq \sqrt{I}$ be ideals in a formally
equidimensional local ring $A$ such that $\ell(I)=\het I$, where $\ell(I)$ denotes the analytic spread of $I$. Then $I$ is
a reduction of $J$ (equivalently $J \subseteq \overline I$) if and only if the $A_\fp$-ideals $I_\fp$ and $J_\fp$ have
the same multiplicity for every minimal prime divisor $\fp$ of
$I$}. 

Using the $j$-multiplicity defined by Achilles and Manaresi
\cite{AM1} (a generalization of the classical Samuel multiplicity),
Flenner and Manaresi \cite{FM} gave  a numerical characterization of
reduction ideals which generalizes B\"{o}ger's result to arbitrary
ideals. 
\begin{flenner}[Flenner-Manaresi \cite{FM}] Let $I \subseteq J$ be
  ideals in a
  formally equidimensional local ring $A$. Then $I$ is a
  reduction of $J$ if and only if $j(I_\fp)=j(J_\fp)$ for all 
$\fp \in \Spec(A)$.
\end{flenner}

It is well known that for an integrally closed domain $A$, the
integral closure of the extended Rees algebra $S=A[It,t^{-1}]$ in
its quotient field is $\overline S=\bigoplus_{n \in \mathbb{Z}}
\overline{I^n}t^n$ ($I^n=A$ for $n<0$), so one could interpret the above
results as numerical
characterizations of the homogeneous components of $\overline S$.  

Our motivation comes from the study of the $S_2$-ification of the same
extended Rees algebra $S=A[It,t^{-1}]$. Under some assumptions on the
ring $A$, $S$ has an $S_2$-ification of the form $\widetilde
S=\bigoplus_{n \in \mathbb{Z}} I_nt^n$, where $I_n =A$ for $n <
0$. In \cite[Theorem 2.4]{C} we proved that if $I$ is primary to
the maximal ideal $\fm$, then $I_n$ is the largest ideal containing
$I^n$ such that $e_i(I_n)=e_i(I^n)$ for $i=0,1$, where $e_0$ and $e_1$
are the first two Hilbert coefficients.

In this paper we use the $j$-multiplicity of Achilles and Manaresi and
a new invariant $j_1$ to obtain a characterization of
$\widetilde S$ similar to the one of  $\overline S$ given by the result
of Flenner and Manaresi.         

The paper is organized as follows. In  the introductory section we
establish the notation and recall the main concepts used in the paper.      

In the second section  we define a generalization of the classical
Hilbert coefficients. Achilles and Manaresi \cite{AM1} defined the so-called
$j$-multiplicity of an ideal $I$ in a local ring $A$ which
generalizes to ideals of maximal
analytic spread the classical Samuel multiplicity. In a subsequent paper,  
  Achilles and Manaresi \cite{AM2} also observed
  that this new invariant can be recovered from the Hilbert polynomial 
  of the bigraded ring $G_\fm(G_I(A))$. 

This is the point of view we adopt in order to define the coefficients 
$j_k(I)\in \mathbb{Z}^{k+1}$ ($0\leq k \leq \dim A$), a
generalization of the classical Hilbert coefficients $e_k(I)$. When
the ideal $I$
is $\fm$-primary, $j_k(I)=(0,\ldots,0,(-1)^ke_k(I))$.  
We show that these coefficients behave
well with respect to general hyperplane sections, one of the main
properties one might expect from any generalization of the Hilbert
coefficients.

The concept of first coefficient ideals has been introduced by Shah in
\cite{Sh}. He proved that for an $\fm$-primary ideal $I$ in a formally
equidimensional ring $(A,\fm)$ there
exists a unique ideal $I_{\{1\}}$, the first coefficient ideal of $I$,
that is maximal among the ideals containing $I$ for which the first
two Hilbert coefficients are equal to those of $I$. In Section 3 we  extend the
definition of
$I_{\{1\}}$ to not necessarily $\fm$-primary ideals. Our definition is
a slight reinterpretation (but necessary for our purpose) of a
description of the first coefficient ideals given by Shah.

  We then observe that using the new
  definition of $I_{\{1\}}$ for an
  arbitrary ideal, 
  we also have $I_n=(I^n)_{\{1\}}$ ($\widetilde S =\bigoplus_{n \in\mathbb{Z}}
  I_nt^n$ is the $S_2$-ification of the extended Rees algebra $S$).  
  This  follows from the proof of \cite[Theorem 2.4]{C} as a  direct
  consequence of an argument due to Heinzer and 
  Lantz \cite[2]{HL}.
  
  The last section contains the main result of this paper. We
  give a numerical characterization of the
  homogeneous components of $\widetilde S$ by proving the following theorem.
\begin{main} Let $\ringA$ be a formally equidimensional local ring and
  let $I \subseteq J$ be ideals of positive height. Then the following are 
   equivalent.
\begin{itemize}
\item[(1)] $J \subseteq I_{\{1\}}$. 
\item[(2)] $j_0(I_\fp)= j_0(J_\fp)$ and $j_1(I_\fp)= j_1(J_\fp)$ for
  all $\fp \in \Spec(A)$.
\end{itemize} 
\end{main}
\noindent Here $j_0(I)=j(I)$ is the above mentioned $j$-multiplicity. 

In fact, we prove a more general version for modules (but technically simpler
for our inductive argument). The proof of the theorem in the
$2$-dimensional case is a crucial part of the argument (see
\ref{fin_lemma}, \ref{lowdim}, and \ref{mainth}). 

%
%
%

\section{Preliminaries}
Throughout this paper a local ring $\ringA$ will be a commutative 
Noetherian ring with identity, and unique maximal ideal. 
\begin{notation}\label{not} 
  Let $\ringA$ be a local ring, let 
  $I$ be an ideal
  of $A$, and let $M$ be a finitely generated $A$-module of dimension $d$. We
  consider the associated graded ring 
  $$G_I(A):=\bigoplus_{n \geq 0}I^n/I^{n+1},$$
  and the associated graded module
  $$G_I(M):=\bigoplus_{n\geq 0}I^nM/I^{n+1}M. $$ 

Given $g \in M\setminus\{0\}$, let $n$ be the largest number such that $g\in I^nM$,
and define the initial form of $g$, denoted $g^{*}$, by
$$ g^{*}:=g\quad\text{modulo}\quad I^{n+1}M \in I^nM/I^{n+1}M
\subseteq  G_I(M).$$ If $g=0$, we define $g^{*}=0$. 
For an $A$-submodule $N$ of $M$,
$$G_I(N,M):=\bigoplus_{n \geq 0} ((N \cap I^nM)+ I^{n+1}M)/I^{n+1}M$$
will denote the $G_I(A)$-submodule of $G_I(M)$ generated by the initial forms of
all elements of $N$.

If the length $\length(M/IM)$ is finite, then for sufficiently large
values of $n$,  $\length (M/I^nM)$ is a polynomial
$P_{I}^{M}(n)$ in $n$ of degree $d$, the Hilbert polynomial of $(I,M)$. 
We write this polynomial in terms of binomial  coefficients:
$$ P_{I}^{M}(n)=e_0(I,M)\binom{n+d-1}{d}-e_1(I,M)\binom{n+d-2}{d-1}+\cdots+(-1)^de_d(I,M).$$ 
The coefficients $e_i(I,M)$ are integers and we call them the Hilbert
coefficients of $(I,M)$.
\end{notation}

\begin{s2}If $A$ is a  Noetherian
  ring, we say that a finitely generated $A$-module $M$ satisfies Serre's
  $(S_2)$ property if for every prime ideal $\fp$ of $A$,
$$\depth M_{\fp} \geq \inf \{2, \dim M_{\fp} \} .$$ 
We say that the ring $A$ satisfies $(S_2)$ if it satisfies   $(S_2)$
as an $A$-module, i.e., $A$ has no embedded prime ideals and 
$\het \fp =1$ for all $\fp \in \Ass(A/xA)$ for any regular element $x \in A$.  
\end{s2}
\noindent We recall the definition of the $S_2$-ification of a Noetherian domain. 
\begin{Definition}
Let $A$ be a Noetherian domain. We say that a domain $B$ is
an $S_2$-ification of $A$ if
\begin{itemize} 
\item[(1)] $A \subseteq B \subseteq Q(A)$ and $B$  is module-finite over $A$,
\item[(2)] $B$ is $(S_2)$ as an $A$-module, and 
\item[(3)] for all $b$ in $B \setminus A$, $\het D(b) \geq 2$, where $
  D(b)=\{a \in A \mid ab \in A\}$.
\end{itemize}
\end{Definition}
\begin{Remark}(\cite[2.4]{HH})
Set $C:= \{b  \in Q(A)\mid  \het D(b) \geq 2 \}$. Then $A$ has an 
$S_2$-ification if and only if $C$ is a
finite extension of $A$, in which case $\tilde A =C$. It is also easy to observe that $\widetilde A$ is 
a finite extension of $A$ inside the quotient field, minimal with
the property that it has the $(S_2)$ property as an $A$-module. 
\end{Remark} 
\begin{Remark}\label{ex_s2} The $S_2$-ification does 
exist for a large class of Noetherian domains. For instance, if $A$ is a universally 
catenary, analytically unramified domain, then
$A$ has an $S_2$-ification
(\cite[EGA,5.11.2]{G}). Also, for any local domain   $\ringA$ that has
a canonical module 
  $\omega$, $A \into \Hom{A}{\omega}{\omega}$ is
  an $S_2$-ification of $A$ (\cite[2.7]{HH}).
\end{Remark} 
We refer to \cite{G}, \cite{A1}, \cite{A2}, and \cite{HH} for  more results about $S_2$-ification.
\begin{coef}
Shah (\cite[Theorem 1]{Sh})
  has proved that if $I$ is an ideal primary to the maximal ideal of a
  formally equidimensional local ring $\ringA$, then 
  the set $$\{J\mid J \text{ ideal of } A, J\supseteq I, e_i(I,A)=
  e_i(J,A) \text
  { for } i=0,1\}$$ has a unique
  maximal element $I_{\{1\}}$, the \emph{first coefficient ideal}
  of $I$. For more about the structure and properties of first
  coefficient ideals we refer the reader to the original paper of Shah
  \cite{Sh} and the series of papers of Heinzer, Lantz, Johnston, and
  Shah (\cite{HJL}, \cite{HJLS}, \cite{HL}).
\end{coef}
In \cite{C} we have proved  the following result:
\begin{Theorem}[\cite{C} Theorem 2.5 and Lemma 2.4]\label{main1}
 Let $\ringA$ be a formally equidimensional,
  analytically unramified local domain with infinite residue field
  and positive dimension, and let $I$ be an $\fm$-primary ideal of $A$. Let
  $\widetilde S=\bigoplus_{n \in \mathbb{Z}}I_nt^n$ be the $S_2$-ification
  of $S=A[It,t^{-1}]$. Then $$I_n \cap A=(I^n)_{\{1\}}\quad \textrm{ for all }n \geq
  1.$$ If $A$ has the $(S_2)$ property, then $I_n$ is an ideal of $A$,
  hence $I_n=(I^n)_{\{1\}}$ for all $n \geq 1$.
\end{Theorem}
\begin{hilb}\label{hilb} 
We first introduce some known facts about Hilbert functions of bigraded
modules. For a detailed description of their properties and  complete
proofs we refer the reader to \cite{D}, \cite{W1}, and \cite{W2} (in
these papers the theory is developed for bigraded rings but it can be
easily extended to bigraded modules).

Let $R=\bigoplus_{i,j=0}^{\infty}R_{ij}$ be a bigraded ring and
let $T=\oplus_{i,j=0}^{\infty}T_{ij}$ be a bigraded $R$-module. Assume
that $R_{00}$ is Artinian  and that $R$ is finitely generated
as an $R_{00}$-algebra by elements of $R_{01}$ and $R_{10}$. 
The Hilbert function of $T$ is defined to be 
$$h_T(i,j)=\length_{R_{00}}(T_{i,j}) .$$
For $i,j$ sufficiently large, the function $h_T(i,j)$ becomes a
polynomial $p_T(i,j)$. If $d$ denotes the dimension of the module  $T$,
we can write this polynomial in the form
$$p_T(i,j)=\sum_{{k,l \ge 0} \atop {k+l \le d-2}} a_{k,l}(T) {i+k \choose k}{j
  +l\choose l},$$ 
with $a_{k,l}(T)$ integers and $a_{k,d-k-2}(T) \ge 0$.

We also consider the sum transform of $h_T$ with respect to the first variable defined by
$$h_T^{(1,0)}(i,j)=\sum_{u=0}^{i}h_T(u,j),$$ and the sum transform of $h_T^{(1,0)}$ with respect to the second variable,
 $$h_T^{(1,1)}(i,j)=\sum_{v=0}^{j}h_T^{(1,0)}(i,v)
=\sum_{v=0}^{j}\sum_{u=0}^{i}h(u,v).$$
For $i,j$ sufficiently large, $h^{(1,0)}(i,j)$ and $h^{(1,1)}(i,j)$
become polynomials with rational coefficients of degrees at most $d-1$
and $d$ respectively. As usual, we can write  these polynomials in
terms of binomial coefficients
$$p_T^{(1,0)}(i,j)=\sum_{{k,l \ge 0}\atop{k+l \le d-1}} a_{k,l}^{(1,0)}(T) {i
  +k\choose k}{j+l \choose l}, $$ 
with $a_{k,l}^{(1,0)}(T)$ integers and $a_{k,d-k-1}^{(1,0)}(T) \geq
0$, and
$$p_T^{(1,1)}(i,j)=\sum_{{k,l \ge 0}\atop{k+l \le d}} a_{k,l}^{(1,1)}(T) {i
  +k\choose k}{j+l \choose l}, $$ 
with $a_{k,l}^{(1,1)}(T)$ integers and $a_{k,d-k}^{(1,1)}(T)\geq 0$.

Since 
$$ h_T(i,j)=h_T^{(1,0)}(i,j)-h_T^{(1,0)}(i-1,j),$$ we get
\begin{equation}\label{dif-eq1}
a_{k+1,l}^{(1,0)}(T)=a_{k,l}(T)\quad \textrm{for } k,l \ge 0, k+l \le
d-2.
\end{equation}

 Similarly we have
$$h_T^{(1,0)}(i,j)=h_T^{(1,1)}(i,j)-h_T^{(1,1)}(i,j-1),$$
 which implies that
\begin{equation}\label{dif-eq2}
a_{k,l+1}^{(1,1)}(T)=a_{k, l}^{(1,0)}(T) \quad \textrm{for } k,l \geq
0, k+l \leq d-1.
\end{equation}
 
\end{hilb}

\section{Generalized Hilbert coefficients}
In this section we define Hilbert coefficients for an arbitrary ideal
$I$ in a local ring $(A,\fm)$. The $k^{\mathrm{th}}$
generalized Hilbert coefficient $j_{k}(I)$ is an element of
$\mathbb{Z}^{k+1}$ whose first $k$ components are $0$ when the ideal
$I$ is primary to the maximal ideal $\fm$. We also show that 
sufficiently general hyperplane sections
behave well with respect to the generalized Hilbert coefficients. This
is one of the main properties that one would expect from a ``good''
definition of these coefficients. 

Let $(A,\fm)$ be a local ring, let $I$ be an ideal of
$A$, and let $M$ be a finitely generated $A$-module of
dimension $d$. Consider the bigraded ring $R=G_{\fm}(G_I(A))$ and the
bigraded $R$-module $T=G_{\fm}(G_I(M))$, where the graded components are
\begin{align*} 
R_{ij}&=(\fm^{i}I^j+I^{j+1})/(\fm^{i+1}I^j+I^{j+1}) \quad \text{and}\\
T_{ij}&=(\fm^{i}I^jM+I^{j+1}M)/(\fm^{i+1}I^jM+I^{j+1}M), \text{
  respectively}.
\end{align*}
Observe that $R_{00}=A/\fm$ and $\dim T=\dim M =d$.

As described in \ref{hilb}, we define  the polynomials
$p_R^{(1,0)}(i,j)$, $p_R^{(1,1)}(i,j)$, $p_T^{(1,0)}(i,j)$, and
$p_T^{(1,0)}(i,j)$. Note that  for $i,j \gg 0$ 
\begin{align*}
p_R^{(1,0)}(i,j)&=\length \big(I^j/(\fm^{i+1}I^{j}+I^{j+1})\big) \quad
\text{and}\\
p_T^{(1,0)}(i,j)&=\length \big(I^jM/(\fm^{i+1}I^{j}M+I^{j+1}M)\big).
\end{align*}
\begin{Definition}\label{coef-def}Let $(A,\fm)$  be a local ring, 
let $I$ be an ideal of  $A$, and let 
    $M$ be a finitely generated $A$-module. Using the notation introduced in \ref{hilb},
    we define 
$$j_k(I,M):=(a_{k,d-k}^{(1,1)}(T),
a_{k-1,d-k}^{(1,1)}(T),\ldots,a_{0,d-k}^{(1,1)}(T)) \in \mathbb{Z}^{k+1} \quad \textrm{for }
0 \leq k \leq d,$$ 
and  call them the {\it generalized Hilbert coefficients} of $(I,M)$.
\end{Definition}
Our main concern will be with the first two coefficients
\begin{align*}
j_0(I,M)&=a_{0,d}^{(1,1)}(T) \quad \textrm{and}\\
j_1(I,M)&=(a_{1,d-1}^{(1,1)}(T), a_{0,d-1}^{(1,1)}(T)).
\end{align*}
To simplify the notation, we denote $j_1(I,M)=(j_1^1(I,M),j_1^2(I,M))$.

\begin{Remark}\label{simple}We also have 
\begin{align*}
j_0(I,M)&=a_{0,d-1}^{(1,0)}(T), \\
j_1(I,M)&=(a_{1,d-2}^{(1,0)}(T), a_{0,d-2}^{(1,0)}(T)), \\
&\ldots\\
j_{d-1}(I,M)&=(a_{d-1,0}^{(1,0)}(T),
a_{d-2,0}^{(1,0)}(T),\ldots,a_{0,0}^{(1,0)}(T)).
\end{align*}
This follows from the equalities (\ref{dif-eq1}) and (\ref{dif-eq2}).
Note that we need to assume $d=\dim M \geq 2$ in order to refer to $j_1(I,M)$ as $(a_{1,d-2}^{(1,0)}(T), a_{0,d-2}^{(1,0)}(T))$. For technical reasons (see Proposition \ref{dif}), we will prefer this interpretation of the generalized Hilbert coefficients. (We only need $d=\dim M \geq 1$ in order to see $j_1(I,M)$ as $(a_{1,d-1}^{(1,1)}(T), a_{0,d-1}^{(1,1)}(T))$.) 

\end{Remark}

\begin{Remark}The coefficients we defined are  a generalization
of the classical Hilbert coefficients. Indeed, when $I$ is $\fm$-primary, 
$$j_k(I,M)=(0,0,\ldots,0,(-1)^ke_k(I,M))\in \mathbb{Z}^{k+1} \quad
\textrm{for } 0 \leq k \leq d, $$ 
  where the first $k$ components are $0$ and  $e_k(I,M)$ is the $k^{\mathrm{th}}$ Hilbert  coefficient of
  $(I,M)$. To see this, note that if  $I$ is
  $\fm$-primary, there exists
  $t$ such  that $\fm^t \subseteq  I$, and then, for $i,j$ large enough,
  $$p_T^{(1,1)}(i,j)=\length(M/I^{j+1}M).$$ An elementary
  identification of the coefficients gives the above equalities.
\end{Remark}

\begin{jmult}\label{jmult} Achilles and Manaresi \cite{AM1} defined
a multiplicity for ideals of maximal analytic spread 
that generalizes
the classical Samuel multiplicity. For a detailed presentation of this 
multiplicity we refer the reader to \cite[Chap. 6]{FCV}.     

Let $(A,\fm)$ be a local
ring, let $I$ be an ideal, and let $M$ be a finitely
generated $A$-module. Then $H_\fm^{0}(G_I(M))$ is a graded
$G_I(A)$-submodule of $G_I(M)$  and is annihilated by $\fm^k$ for $k$
large enough, so it may be considered as a module over
$\bar G_I(A) := G_I(A) \otimes_A A/{\fm^k}$. Then $e(H_\fm^{0}(G_I(M)))
:=e(\bar G_I(A)^{+},H_\fm^{0}(G_I(M)))$ is well defined, where
$\bar G_I(A)^{+}$ denotes the ideal of $\bar G_I(A)$ of elements of
positive degree. Thus we can define
\begin{displaymath}
j(I,M):=\left\{ \begin{array}{ll}
     e(H_\fm^{0}(G_I(M)) & \textrm{if } \dim M = \dim H_\fm^{0}(G_I(M))\\
     0 & \textrm{if } \dim M >  \dim H_\fm^{0}(G_I(M))
        \end{array} \right. 
\end{displaymath} 
\noindent 
Note that $j(I,M)\not=0$ if and only if $\ell_M(I)=\dim M$
\cite[6.1.6(1)]{FCV},
where $\ell_M(I)=\dim G_I(M)/\fm G_I(M)$ (the analytic spread of $I$
in $M$). 
\end{jmult} 
\begin{genm}
In \cite{AM2} Achilles and Manaresi defined another generalization of
the Samuel multiplicity. Our presentation will be given
in the slightly more general context of modules. 

Let $I$ be an arbitrary ideal in a local ring $(A,
\fm)$, and let $M$ be a finitely generated $A$-module.
Using the notation introduced in \ref{hilb}, denote 
$$c_i(I,M):=a_{i,d-i}^{(1,1)}(T) \quad (0 \leq i \leq d),$$
where $T=G_\fm(G_I(M))$. The sequence $(c_i(I,M))_{0 \leq i \leq d}$ is
called the multiplicity sequence of $(I,M)$. In the case $M=A$ we
simply denote $c_i=c_i(I,A)$.

Note that this sequence consists of the leading coefficients of the
generalized Hilbert coefficients that we defined in \ref{hilb}.

We state the  following proposition proved in \cite{AM2} (we present a
version for modules).  
\begin{Proposition}[{\rm{\cite[Proposition
      2.3]{AM2}}}]\label{am-prop}  
Let $(A,\fm)$ be a local ring, let  $I$ be a proper
ideal of $A$, and let $M$ be a finitely generated $A$-module. 
Set $l=\dim G_I(M)/\fm G_I(M)$  and $q=\dim(M/IM)$. Then
\begin{itemize}
\item[\rm(i)] $c_k(I,M)=0$ for $k < d - l$ or $ k > q$;
\item[\rm(ii)] $c_{d-l}(I,M)=\sum_{\mathfrak{\beta}}e(\fm
G_{\mathfrak{\beta}},G_I(M)_\beta)e(G/\mathfrak{\beta})$,
where $\mathfrak{\beta}$ runs through the all highest dimensional
associated primes of $G_I(M)/{\fm G_I(M)}$ such that \newline $\dim(G/\mathfrak{\beta})+
\dim G_{\mathfrak{\beta}}=\dim G$;
\item[\rm(iii)] $c_q(I,M)=\sum_{\fp}e(IA_\fp,M_{\fp})e(A/\fp),$
where $\fp$ runs through the all highest dimensional
associated primes of $M/IM$ such that $\dim A/\fp+ \dim A_\fp=\dim A.$
\end{itemize}
\end{Proposition}

Achilles and Manaresi \cite[Proposition
2.4]{AM2} also proved that the $j$-multiplicity $j(I,M)$ is equal to the 
coefficient $c_0(I,M)$.
For more details   we refer the reader to the original paper of
Achilles and Manaresi \cite{AM2} (the proofs can be immediately
extended to the version for modules we present here). 

We will prove  that the multiplicity sequence defined
above is an invariant of the ideal up to its integral closure. If $J
\subseteq I$, we say that $J$ is a reduction of $(I,M)$ if there
exists $n$ such
that $JI^nM=I^{n+1}M$.
\begin{Proposition}\label{lat} Let $(A,\fm)$ be a local ring, let $J \subseteq I$
  be proper  ideals
  of $A$, and let $M$ be a finitely generated $A$-module. If $J$ is a
  reduction of $(I,M)$, then $c_i(J,M)=c_i(I,M)$ for $i=0,\ldots,d$.
\end{Proposition}

Since the proof requires technical results that will be made clear later,
we postpone it  until the end of this paper.
\end{genm}

Before proceeding further, we need to introduce more notation.

If $x$ is an element of  $A$, denote by $x'$ the initial form
of $x^{*} \in G_I(A)$ in $R=G_\fm(G_I(A))$. Similarly, if $J$ is an
ideal in $A$, let $$J'= G_\fm(G_I(J,A),G_I(A)) \subseteq R$$ 
be the ideal generated
by all $x'$ when $x \in J$, and if $N$ is an $A$-submodule of $M$, we 
denote $$N'=G_\fm(G_I(N,A),G_I(M)) \subseteq T=G_\fm(G_I(M)).$$

\begin{Definition}[\cite{D}]
Let $R=G_\fm(G_I(A))$ and let  
$(0)=N_1 \cap N_2 \cap \ldots \cap N_r \cap N_{r+1}\cap \ldots \cap N_t$ 
be an irredundant primary decomposition of $(0)$ in the $R$-module $T=G_\fm(G_I(M))$. Denote
$P_i=\sqrt{ (N_i:_R T)}$, $i=1,\ldots,t$. Assume that
\begin{align}
\label{eq:pr1}
 I' &\subseteq P_{r+1},\ldots,P_t \quad \textrm{and}\\
\label{eq:pr2}
I' &\not\subseteq P_1, \ldots,P_r.
\end{align}
We say that $x \in I$ is a superficial element for $(I,M)$ if $x' \not\in P_1, \ldots, P_r$. 
\end{Definition}
Note that we can always choose $x \in I \setminus \fm I$ superficial element for $(I,M)$.
\begin{Remark}\label{sup_rem}
Let $x\in I$ be a superficial element for $(I,M)$. By (\ref{eq:pr1}), there exists $k$ such that $({I'})^kT \subseteq N_{r+1} \cap \ldots \cap
N_t$. Then  
$$(0':_T x')=\bigcap_{i=1}^{t} (N_i:_T x') \subseteq N_1 \cap \ldots
\cap N_r, $$
hence
\begin{equation}
\label{eq:pr3}
({I'})^kT \cap (0':_{T}x') \subseteq  N_1 \cap N_2 \cap \ldots \cap N_r \cap N_{r+1}\cap \ldots \cap N_t=(0).
\end{equation}
\end{Remark}
The following lemma, in its  version for ideals, is due to Dade
\cite[3.1]{D}(unpublished thesis). For convenience, we present here a proof. 
\begin{Lemma}\label{dade_lemma} Let $A$ be a Noetherian ring, let $I$ be an
  ideal of $A$, let $M$
  be a finitely generated $A$-module, and
  let $L \subseteq K$ be two submodules of $M$ such that the length $\length
  (K/L)$ is  finite. Then $$ \length (K/L)= \length
  (G_I(K,M)/G_I(L,M)).$$ 
\end{Lemma}
\begin{proof} 
Consider the descending chain of modules
\begin{displaymath}
\frac{K \cap IM +L}{L} \supseteq \frac{K \cap I^2M +L}{L} \supseteq
\ldots .
\end{displaymath}
The module $K/L$ has finite length, so there exists $N$ such that 
$$\frac{K \cap I^nM +L}{L}=\frac{K \cap I^{n+1}M +L}{L} \quad
\text{for } n > N $$
which implies that 
$$ K \cap I^nM +L=K \cap I^{n+1}M + L \quad \text{for } n > N. $$
So, for $n > N$,
$$K \cap I^nM +L \subseteq \bigcap_{k \geq 1}(K \cap I^kM +L)
\subseteq \bigcap_{k \geq 1}(I^kM+L)=L,$$
i.e.,  $K \cap I^nM = L  \cap I^nM$.

Finally,
\begin{align*}
\length(K/L)&=\length\Big(\frac{K+IM}{L+IM}\Big)+\length\Big(\frac{K
  \cap IM}{L\cap IM}\Big)\\
            &=\length\Big(\frac{K+IM}{L+IM}\Big)+ \length\Big(\frac{K
  \cap IM+I^2M}{L\cap IM+I^2M}\Big) + \length \Big(\frac{K
  \cap I^2M}{L\cap I^2M}\Big)\\
&\cdots\\
&=\length\Big(\frac{K+IM}{L+IM}\Big)+\cdots+
  \length\Big(\frac{K+I^NM}{L+I^NM}\Big)\\
&=\length\Big(\frac{G_I(K,M)}{G_I(L,M)}\Big).
\end{align*}\end{proof}

The following proposition shows that sufficiently general hyperplane
sections behave well with respect to the generalized Hilbert coefficients. 
\begin{Proposition}
\label{dif}
Let $(A,\fm)$ be a local ring and let $M$ be a finitely generated
$A$-module. Suppose that $x \in I$ is a superficial element for $(I,M)$  and a
nonzerodivisor on $M$ with $x'
\in R_{01}$. Denote $\overline T= G_{\overline \fm}(G_{\overline I}(\overline M))$, where $\overline
A=A/xA$, $\overline I= I \otimes_A \overline A$, and $\overline M = M
\otimes_A \overline A$.
Then, for $i,j$ large,  
$$h_T^{(1,0)}(i,j)-h_T^{(1,0)}(i,j-1)= h_{\overline T}^{(1,0)}(i,j) .$$
In particular, $j_0(I,M)=j_0(\overline I,\overline{M}), j_1(I,M)=j_1(\overline I,\overline
M),\ldots,j_{d-1}(I,M)=j_{d-1}(\overline I,\overline M)$, where $d$ denotes the
dimension of the module $M$.
\end{Proposition}

\begin{proof}
The proof relies on  Lemma \ref{dade_lemma}, a technique also used by
Dade in \cite{D}.

We have the following exact sequence
$$0 \to K \to \frac{I^j M}{\fm^{i+1}I^j M+I^{j+1} M} \to
\frac{I^jM+xM}{\fm^{i+1}I^jM+I^{j+1}M+xM} \to 0,$$
where 
\begin{align*}
K&=\frac{I^j M \cap (\fm^{i+1}I^jM+I^{j+1}M+xM)}{\fm^{i+1}I^j
M+I^{j+1} M}\\&= \frac{(\fm^{i+1}I^jM+I^{j+1}M)+ I^j M \cap xM}{\fm^{i+1}I^j
M+I^{j+1} M}\\
&\cong\frac{I^j M \cap xM}{(\fm^{i+1}I^jM+I^{j+1}M)\cap xM} \quad .            \end{align*}

From this exact sequence we get 
\begin{align*}
h_{\overline T}^{(1,0)}(i,j)&= \length \Big(\frac{I^jM+xM}{\fm^{i+1}I^jM+I^{j+1}M+xM}\Big)\\
&= \length \Big(\frac{I^j M}{\fm^{i+1}I^j M+I^{j+1} M}\Big) - \length
\Big(\frac{I^j M \cap xM}{(\fm^{i+1}I^jM+I^{j+1}M)\cap xM}\Big).
\end{align*}
 Therefore we need to prove that for $i,j \gg 0$
\begin{displaymath}
\length \Big(\frac{I^j M \cap xM}{(\fm^{i+1}I^jM+I^{j+1}M)\cap xM}\Big)=\length
  \Big(\frac{I^{j-1}M}{\fm^{i+1}I^{j-1}M+I^jM} \Big).
\end{displaymath}
We have
\begin{align*}
\length \Big(\frac{I^j M\cap xM}{(\fm^{i+1}I^jM+I^{j+1}M)\cap
  xM}\Big)&=\length
  \Big(\frac{x(I^jM:x)}{x((\fm^{i+1}I^jM+I^{j+1}M):x)}\Big)\\
&=\length  \Big(\frac{(I^jM:x)}{(\fm^{i+1}I^jM+I^{j+1}M):x}\Big)\\
&=\length
  \Big(\frac{(I^jM:x)'}{((\fm^{i+1}I^jM+I^{j+1}M):x)'}\Big),
\end{align*}
where the last equality follows by a successive application of Lemma \ref{dade_lemma}.

By Remark \ref{sup_rem}, there exists  $c$ such that  $({I'})^cT \cap (0':_{T}x')=(0)$.
We claim that for $j > c$
\begin{align}
\label{cond1}
(I^jM:x)'\cap ({I'})^{c}T&=(I^{j-1}M)' \quad 
\text{ and}\\
\label{cond2}
((\fm^{i+1}I^jM+I^{j+1}M):x)' \cap
({I'})^{c}T&=(\fm^{i+1}I^{j-1}M+I^jM)' .
\end{align} 

We first prove (\ref{cond1}). Let $y \in (I^jM:x)$ such that $0 \neq y'
\in ({I'})^cT$. Since $({I'})^cT \cap (0':_{T}x')=(0)$, it follows
that $y' \notin (0':x')$, hence $0 \neq (yx)' \in (I^jM)'$. But
$(I^jM)'$ is 
\[ \begin{array}{*{13}c}
0 &\oplus& 0 & \oplus & \cdots & \oplus & 0 & \oplus &
T_{0,j} & \oplus & T_{0,j+1} & \oplus & \cdots \\
\oplus& & \oplus & &  & & \oplus & & \oplus & & \oplus & &\\ 
0 &\oplus& 0 & \oplus & \cdots & \oplus & 0 & \oplus
& T_{1,j} & \oplus & T_{1,j+1} & \oplus & \cdots \\
\oplus& & \oplus & &  & & \oplus & & \oplus & & \oplus & &\\ 
0 &\oplus& 0 & \oplus & \cdots &\oplus & 0 & \oplus &
T_{2,j} & \oplus & T_{2,j+1} & \oplus & \cdots \\
\oplus& & \oplus & &  & & \oplus & & \oplus & & \oplus & &\\ 
\vdots& & \vdots & &  & & \vdots & & \vdots & & \vdots & &
\end{array} \]
Since $x'
\in R_{01}$, we must have $y' \in (I^{j-1}M)'$.

To see (\ref{cond2}), consider $y \in ((\fm^{i+1}I^jM+I^{j+1}M):x)$
such that $0 \neq y' \in ({I'})^cT$. By the choice of $c$, we have 
$y' \notin (0':x')$, hence $(yx)' \in ( \fm^{i+1}I^jM+I^{j+1}M)'$ and $(yx)' \neq 0$. The
homogeneous components of the graded
submodule $(\fm^{i+1}I^jM+I^{j+1}M)' \subseteq T$ are represented
below:
\[ \begin{array}{*{13}c}
0 & \oplus & \cdots & \oplus & 0 & \oplus & 0 & \oplus &
T_{0,j+1}
& \oplus & T_{0,j+2} &  \oplus & \cdots\\
\oplus & &  & & \oplus & & \oplus & & \oplus & & \oplus& &  \\
0 & \oplus & \cdots & \oplus & 0 & \oplus & 0 & \oplus &
T_{1,j+1}
& \oplus & T_{1,j+2} & \oplus & \cdots\\
\oplus & &  & & \oplus & & \oplus & & \oplus & & \oplus& & \\
\vdots & &  & & \vdots & & \vdots & & \vdots & & \vdots& &  \\
\oplus & &  & & \oplus & & \oplus & & \oplus & & \oplus& &   \\
0 & \oplus & \cdots & \oplus & 0 & \oplus & 0 & \oplus &
 T_{i,j+1}
& \oplus & T_{i,j+2} & \oplus &  \cdots\\
\oplus & &  & & \oplus & & \oplus & & \oplus & & \oplus& &  \\
0 & \oplus &\cdots & \oplus &0 & \oplus &T_{i+1,j} & \oplus &T_{i+1,j+1}
& \oplus &T_{i+1,j+2} & \oplus & \cdots\\
\oplus & &  & & \oplus & & \oplus & & \oplus & & \oplus& &
   \\
0 & \oplus &\cdots & \oplus &0 & \oplus &T_{i+2,j} & \oplus &T_{i+2,j+1}
& \oplus &T_{i+2,j+2} &  \oplus & \cdots\\
\oplus & &  & & \oplus & & \oplus & & \oplus & & \oplus& &  \\
\vdots & &  & & \vdots & & \vdots & & \vdots & & \vdots& &
\end{array} \]
Since $x' \in R_{01}$ we get $y' \in
(\fm^{i+1}I^{j-1}M+I^jM)'$. 

Then we have
\begin{align*}
&\length \Big(\frac{I^jM \cap xM}{(\fm^{i+1}I^jM+I^{j+1}M)\cap
  xM}\Big)\\
&= \length
\Big(\frac{x(I^jM:x)}{x((\fm^{i+1}I^jM+I^{j+1}M):x)}\Big)\\
&=\length
\Big(\frac{(I^jM:x)}{((\fm^{i+1}I^jM+I^{j+1}M):x)}\Big)\\
&=\length
\Big(\frac{(I^jM:x)'}{((\fm^{i+1}I^jM+I^{j+1}M):x)'}\Big)\\
&=\length
\Big(\frac{(I^jM:x)'\cap ({I'})^cT}{((\fm^iI^jM+I^{j+1}M):x)'\cap ({I'})^cT}\Big)+
 \length \Big(\frac{(I^jM:x)' + ({I'})^cT}{((\fm^{i+1}I^jM+I^{j+1}M):x)' +
  ({I'})^cT}\Big)\\
&=\length \Big(\frac{I^{j-1}M }{(\fm^{i+1}I^{j-1}M+I^{j}M)}\Big)+  \length \Big(\frac{(I^jM:x)' + ({I'})^cT}{((\fm^{i+1}I^jM+I^{j+1}M):x)' +
  ({I'})^cT}\Big).
\end{align*}
By the Artin-Rees lemma, there exists $p$ such that for $j > p$
$$I^jM \cap xM =I^{j-p}(I^pM \cap xM),$$ 
i.e.,
$$x(I^jM :_{M} x)=xI^{j-p}(I^pM :_{M}  x), $$
or
$$(I^jM :_{M}  x)=I^{j-p}(I^pM :_{M} x).$$
Then, for $j > p+c$, 
$(I^jM:x)' \subseteq {I'}^c T.$ 

On the other hand, we also have
\begin{align*}
((\fm^{i+1} I^jM+I^{j+1}M):x)' &\subseteq (I^jM:x)' \\&\subseteq ({I'})^c T \quad \quad \text{for } j >
n+c \quad \text{and all } i.
\end{align*}
We can now conclude that 
  \begin{displaymath}
  \length \Big(\frac{I^j M \cap xM}{(\fm^{i+1}I^jM+I^{j+1}M)\cap
  xM}\Big)=\length
  \Big(\frac{I^{j-1}M}{\fm^{i+1}I^{j-1}M+I^jM} \Big),
  \end{displaymath}
  which finishes the proof.
\end{proof}

\section{First coefficient ideals-the general case}
In this section we define the first coefficient ideal $I_{\{1\}}$ of a
not necessarily $\fm$-primary ideal $I$. We then observe
that using the new definition of $I_{\{1\}}$, Theorem \ref{main1} is
true in general, without assuming that $I$ is
$\fm$-primary. 

For reasons that will become obvious later, we need
again to introduce the notion in the more general context of modules.
 
\begin{Definition}\label{def-coef}Let $M$ be a finitely generated $A$
  module and 
  let $I$ be an ideal of  $A$ with $\dim M/IM < \dim M$. We define
  $I_{\{1\}}^M$, the first 
   coefficient ideal of $(I,M)$, to be the ideal of $A$
   $$I_{\{1\}}^M=  \bigcup (I^{n+1}M :_{A} aM), $$ where the union ranges
   over 
   all $n \ge  1$ and  all $a \in I^n \setminus I^{n+1}$
   such that $a^{*}$ is part of a system of parameters
   of $G_I(M)$. If $M=A$, we simply denote $I_{\{1\}}^M=I_{\{1\}}$. 
\end{Definition}
\begin{Remark} Let us observe that our definition coincides with the one 
given by
  Shah in the $\fm$-primary case. Indeed, by the structure theorem for
  the coefficient ideals proved by Shah (\cite[Theorem 2]{Sh}), we
  have
\begin{equation}\label{shah-descr}
I_{\{1\}}=\bigcup (I^{n+1} :_{A} a), 
\end{equation}
where the union ranges
   over 
   all $n \ge  1$ and  all $a$ extendable to some minimal reduction of
   $I^n$. 

On the other hand, $a$ is  extendable to some minimal reduction of
   $I^n$ if and only if the image of $a^{*}$ in $G_I(A)/\fm G_I(A)$ is
   part of a system of
   parameters. But if the ideal $I$ is $\fm$-primary  this is
   equivalent to the fact that  $a^{*}$ is
   part of a system of
   parameters of $G_I(A)$, for  the ideal $\fm G_I(A)$ is nilpotent.
\end{Remark}
\medskip
Heinzer, Johnston, Lantz, and Shah \cite[Theorem
   3.17]{HJLS}
gave a description of the coefficient ideals
   involving the  blow-up of $I$. We present here their result for
   the case of the first coefficient ideals. 

The blow-up $\mathcal{B}(I)$ of an ideal $I$ in a local 
domain $A$ is defined to be the model
$$\mathcal{B}(I)=\{A[I/x]_{\fp} \mid 0 \neq x \in I \textrm{ and } \fp
\in \Spec (A[I/x])\}.$$ 
$\mathcal{B}(I)$ is the set of all local rings between $A$ and the
quotient field $Q(A)$ minimal with respect to domination among those
in which the extension of $I$
is a principal ideal. Let $\mathcal{D}_{1}$ denote the intersection of
the local domains on the blow-up $\mathcal{B}(I)$  of dimension at
most $1$ in which the maximal ideal is minimal over the extension of
$I$ (see \cite[Definition 3.2]{HJL}). The main result of
\cite{HJLS}(Theorem 3.17) says that if $A$ is a formally
equidimensional, analytically unramified local domain with infinite
residue field and $\dim A > 0$, and $I$ is an $\fm$-primary ideal,
then
\begin{equation}\label{hjls-descr}
I_{\{1\}}=I\mathcal{D}_{1}\cap A.
\end{equation}  

In a subsequent paper, Heinzer and Lantz \cite[2]{HL} prove directly
the equivalence of  the description of the first coefficient ideals
given initially by Shah (see \ref{shah-descr}) and the description
given by \ref{hjls-descr}. The argument assumes that the
ideal $I$ is $\fm$-primary, but a careful examination of their proof
actually shows the following:

\begin{Proposition}
Let $(A,\fm)$ be a formally equidimensional local ring of positive
dimension, and let $I$ be an arbitrary
ideal of $A$. Then 
$$I\mathcal{D}_{1}\cap A=\bigcup (I^{n+1} :_{A} a),$$ 
where the union ranges
   over 
   all $n \ge  1$ and  all $a \in I^n \setminus I^{n+1}$ such that
   $a^{*}$ is part of a system of parameters of $G_I(A)$.
\end{Proposition}
Note that the right hand side of this equality is exactly the
definition of the first
coefficient ideals in the general case (see Definition
\ref{def-coef}).

In \cite{C} we have proved Theorem \ref{main1}. The statement of the
theorem assumes that $I$ is an $\fm$-primary ideal, but all is used in
the proof is that  $I_{\{1\}}=I\mathcal{D}_{1}\cap A$. Therefore, by the above
discussion, we have the following theorem.
\begin{Theorem}\label{main2}
 Let $\ringA$ be a formally equidimensional,
  analytically unramified local domain with infinite residue field
  and positive dimension, and let $I$ be an arbitrary ideal of $A$. If
  $\widetilde S=\bigoplus_{n \in \mathbb{Z}}I_nt^n$ is the $S_2$-ification
  of $S=A[It,t^{-1}]$, then $$I_n \cap A=(I^n)_{\{1\}}\quad \textrm
  { for all }n \geq1, $$ 
where for an ideal $J$, $J_{\{1\}}$ denotes the first coefficient
ideal of $J$ as defined in \ref{def-coef}. 

In particular, if $A$ has the $(S_2)$ property, then $I_n=(I^n)_{\{1\}}$
for all $n \geq 1$.
\end{Theorem}

In this way, the problem of giving a numerical characterization of the
$S_2$-ification of the extended Rees algebra $S=A[It,t^{-1}]$ reduces
to the problem of finding a numerical characterization of the
generalized first coefficient ideals (Definition \ref{def-coef}).

The following proposition shows that the union involved in Definition
\ref{def-coef} can be replaced by a single colon ideal. It is the
analogue of Theorem 3 of \cite{Sh}. 

Recall that a finitely generated
module $M$ over a local ring $A$ is called equidimensional if
for every minimal prime ideal $\fp$ of $M$ the module $M/\fp M$
has dimension $\dim M$. We also say that $M$ is formally
equidimensional if $\widehat M$ (the completion of $M$ in the
$\fm$-adic topology) is equidimensional as an $\widehat A$-module. If
the ring $A$ is complete and $M$ is equidimensional, then $G_I(M)$ is
also equidimensional (see \cite[18.24]{HIO} and \cite[4.5.6]{BH}). 
\begin{Proposition} \label{sh_lem1}Let $M$ be a finitely generated
  formally equidimensional $A$-module and let $I$ be an ideal of  $A$ 
  such that  $\dim M/IM < \dim M$. Then
  there exist a fixed integer
  $m$ and a  fixed element $x$ of $I^m \setminus I^{m+1}$ with
  $x^{*}$ part of system of parameters of $G_I(M)$ such that 
  $$I_{\{1\}}^M=  (I^{m+1}M :_{A} xM).$$
\end{Proposition} 
\begin{proof}
We can assume that $A$ is complete and that $M$ is
equidimensional.  
Let $N$ be  the $G_I(A)$-submodule of $G_I(M)$ generated by  $I_{\{1\}}^M
M/IM$. By definition, each generator of $N$ is annihilated
by a homogeneous element
of $G_I(A)$ which is part of a system of parameters of $G_I(M)$. By prime
avoidance, we can find a homogeneous element $x^{*}\in I^m/I^{m+1}$
($x \in I^m$) 
that annihilates the entire submodule $N$ 
and which avoids all the minimal primes in the support of
$G_I(M)$. The observation that $G_I(M)$ is equidimensional (implied by
the hypothesis) concludes the proof.  
\end{proof}
\begin{Proposition}\label{l_dim}Let  $M$ be a formally equidimensional
  $A$-module, and  let $I \subseteq J$ be 
  ideals of $A$ such that $\dim M/IM < \dim M$. Then  $J \subseteq
  I_{\{1\}}^M$ if and
  only if $$\dim \bigoplus_{n \geq 0} JI^nM/I^{n+1}M < \dim G_I(M)=\dim
  M .$$ 
\end{Proposition}
\begin{proof}
Indeed, if we denote $L= \bigoplus_{n \geq 0}
  JI^nM/I^{n+1}M $, then, by
  Proposition \ref{sh_lem1}, it follows that $L$
  is annihilated by an element which is part of a system of
  parameters of $G_I(M)$.    
\end{proof}

\begin{Remark} If $M$ is faithful (i.e. $\Ann M=0$) and  $J$ is
    a (minimal) reduction of $(I,M)$, then $J$ is a (minimal) reduction of
    $I$. Indeed, if $I^{n+1}M=JI^nM$ for some $n$, then, by the
    determinant trick, $J$ and $I$ have the same integral
    closure, i.e., $J$ is a reduction of $I$.
\end{Remark} 

In the $\fm$-primary case it is obvious that the ideal $I$ is a
reduction of its first coefficient ideal (by definition). This is
still true in the general case, as the following proposition shows.

\begin{Proposition}\label{red} Let $\ringA$ be a  local
  ring, let $M$ be a finitely generated formally equidimensional $A$-module,
  and let $I$ be an ideal of $A$ such that $\dim M/IM < \dim M$. If $ I\subseteq J
  \subseteq  I_{\{1\}}^M$, then $I$ is a reduction of $(J,M)$.
\end{Proposition}
\begin{proof} 
As usual, we may assume that $A$ is a complete local ring.
First we prove the proposition in the case when $M$ is faithful.
Note that in this case both $A$ and $M$ will be equidimensional,
therefore both $G_I(A)$ and $G_I(M)$ are equidimensional of dimension equal to
   $\dim
  A=\dim M$ (this is implicitly proved in Theorem 4.5.6 of
  \cite{BH}).

Let us observe that for a faithful $A$-module $M$, $\Ann G_I(M)$ is
  a nilpotent ideal of $G_I(A)$.  Indeed, if  $\bar{x} \in I^n/I^{n+1}$ is  an
  element of $G_I(A)$ that
  annihilates $G_I(M)$, then $xM \subseteq I^{n+1}M$, which
  by the determinant trick implies that $x \in \overline{I^{n+1}}$ (here
  $\overline{J}$ denotes the integral closure of the ideal $J$) . If we
  write the equation of integral dependence we get
  $$ x^k+a_1x^{k-1}+...+a_k=0,$$ with $ a_i \in I^{(n+1)i}$. Thus
  $x^{k}= -(a_1x^{k-1}+...+a_k) \in I^{kn+1}$, which implies
  that $\bar{x} \in G_I(A)$ is nilpotent.

By Proposition \ref{sh_lem1}, there exist a fixed integer
  $m$ and a fixed element $a \in I^m \setminus I^{m+1}$ with $a^{*}
  \in G_I(A)$
  part of a system of parameters
  of $G_I(M)$ such that
  $I_ {\{1\} }^M=(I^{m+1}M:_{A}aM)$. Let  $y \in
  (I^{m+1}M:_{A}aM)$. Then $yaM \subseteq I^{m+1}M$, and 
  using the determinant trick we get  
\begin{equation}\label{p1}
ya \in \overline{I^{m+1}}.
\end{equation}
Since $\Ann (G_I(M))$ is nilpotent and $G_I(A)$ is equidimensional,
  $a^{*}$  is  part of
  a system of parameters of
  $G_I(A)$, i.e. $at^m \in S=A[It,t^{-1}]$ is not contained in any
  minimal prime divisor of $t^{-1}S$. 

We claim that from the above assertion
   and (\ref{p1}) it follows that $y \in \overline{I}$. To
  prove this, note that we may also assume that $A$ is a reduced ring.
  Let $\overline
  T=\bigoplus_{n\geq
  0}\overline{I^n\overline{A}}t^n$ be
  the integral closure of $T$ in its total quotient ring. Since the
  ring $\overline A$ is equidimensional (it is a local catenary ring
  satisfying the $(S_2)$ property; see \cite[5.10.9]{G}), the ring
  $\overline{T}/t^{-1}\overline T$ is also equidimensional (implicitly
  proved  in Theorem 4.5.6 of \cite{BH}; note that
  $(\overline{I^n})_{n \geq 0}$ is a Noetherian filtration) and 
  is a finite extension of $T/t^{-1}T$. In particular, 
  any minimal prime of $t^{-1}\overline{T}$ contracts
  back to a minimal prime of $t^{-1}T$. Thus the image of $at^m$ does
  not belong to any associate prime of $t^{-1}\overline{T}$, hence
  $a^{*}$ is a nonzerodivisor on $\overline{T}/t^{-1}\overline{T}$.
  By (\ref{p1}) we get  $y \in \overline{I \overline{A}} \cap
  A=\overline{I}$.
\end{proof}
\section{The main result}

We now prove  two propositions  that will be the main tools for
proving Theorem~\ref{mainth} in dimension 2.
\begin{Proposition} \label{fin_lemma}Let $M$ be a finitely generated
  formally equidimensional 
  $A$-module  of
  dimension $2$, and let  $I \subseteq J$
  be two ideals of $A$ such that $\dim M/IM < \dim M$. If  $J
  \subseteq I_{\{1\}}^M$, then there exist  positive integers $k$ and
  $l$ 
  such that $$\fm^k I^jM \subseteq J^jM \quad \text{for }j\geq l.$$
In particular,
$$\length(J^jM/I^jM)  < \infty \quad \text{for } j\gg0.$$ 
\end{Proposition}
\begin{proof} 
  Denote by $N$ the $G=G_I(A)$-submodule of $G_I(M)$
  generated in degree $0$ by $JM/IM$, i.e.
$$N=\bigoplus_{n\geq 0}JI^nM/I^{n+1}M.$$

 By Proposition \ref{l_dim}, we have  $\dim_G(N) \leq \dim G(M) -
  1=1$, which implies  that $\dim
  _{G_\fm(G)}G_\fm(N) \leq 1$. Since  
  $$G_\fm(N)=\bigoplus_{i,j \geq 0}  \frac{\fm^iJI^jM +
    I^{j+1}M}{\fm^{i+1}JI^jM + I^{j+1}M} \ ,$$ 
  It follows that for $i,j \gg 0$
  $$\length \Big(\frac{\fm^iJI^jM + I^{j+1}M}{\fm^{i+1}JI^jM +
    I^{j+1}M}\Big)$$ is a polynomial of degree $\leq \dim G_\fm(N) -2 \leq -1$,
  so there exist $i_0,j_0$ such that 
$$\length \Big(\frac{\fm^iJI^jM + I^{j+1}M}{\fm^{i+1}JI^jM +
    I^{j+1}M}\Big)=0
  \quad \text{for } i \geq i_0, j \geq j_0.$$
  By Nakayama's lemma we then obtain 
  \begin{equation} \label{eq:fl}\fm^i JI^jM \subseteq I^{j+1}M \quad
    \text{for } i \geq i_0, j \geq j_0.
  \end{equation}
  Since $I$ is a reduction of $(J,M)$ (\ref{red}) there exists $n$ such that
  $I^jJ^nM=J^{n+j}M$ for $j \geq 1$. By (\ref{eq:fl}) it follows that
  $$\fm^{ni}J^nI^jM \subseteq I^{n+j}M \quad \text{for } i \geq i_0, j
  \geq j_0,$$ which in conjunction with the
  previous equality implies that $$\fm^{ni}J^{n+j}M \subseteq I^{n+j}M
  \quad \text{for } i \geq i_0, j \geq j_0.$$ 
Take $k=i_0$ and $l=n+j_0$. 
\end{proof}

\begin{Proposition} \label{lowdim}Let $(A,\fm)$ be a local ring and
  let $M$ be a finitely generated formally equidimensional
  $A$-module of
  dimension $\leq 2$. Consider $I \subseteq J$ two ideals in $A$ with
  $\dim M/IM < \dim M$ such that
  $I \subseteq J \subseteq I_{\{1\}}^M$. Then, for $i,j$ large enough,
\begin{itemize}
\item[\rm 1)] $\length (J^jM/I^jM)$ is a constant;
\item[\rm 2)] $\length (\fm^{i+1}J^jM+J^{j+1}M/\fm^{i+1}I^jM+I^{j+1}M)$ is
  a constant;
\item[\rm 3)] $\length(I^jM/\fm^{i+1}I^jM+I^{j+1}M)=\length(J^jM/\fm^{i+1}J^jM+J^{j+1}M)$.
\end{itemize}
\end{Proposition}
\begin{proof} 
 By  Proposition \ref{fin_lemma}, $\length (J^jM/I^jM)$ is finite for $j$ large
  enough, so for the first part of the proposition we
  can use an
  argument similar (but in module version) to the one used by 
  Shah in the proof of Theorem 2
  of \cite{Sh}. 

Since $I \subseteq J \subseteq I_{\{1\}}^{M}$, $I$ is a reduction of
$(J,M)$ (see Proposition \ref{red}), hence there exists an integer $s$
such that  $I^nJ^sM=J^{n+s}M$ for all $n$. Then we have
\begin{align*}
\length\big(J^{s+n}M/I^{s+n}M\big)&=\length\big(J^{s}I^{n}M/I^{s+n}M\big)\\
 &=\sum_{i=1}^{s}\length\big(J^iI^{n+s-i}M/J^{i-1}I^{n+s-i+1}M\big)\\
 &=\sum_{i=1}^{s}\length\big(J^{i-1}I^{s-i}JI^nM/J^{i-1}I^{s-i}I^{n+1}M\big)\\
 &\leq \sum_{i=1}^{s}c_i \length\big(JI^nM/I^{n+1}M\big)
\end{align*}                                      
where $c_i$ is the number of generators of $J^{i-1}I^{s-i}M$. Set
$c=\sum c_i$. Then
$$ \length\big(J^{s+n}M/I^{s+n}M\big) \leq c \length\big(JI^nM/I^{n+1}M\big).$$

On the other hand, by Proposition \ref{l_dim}, for $n$ large enough,
$\length(JI^nM/I^{n+1}M)$ is a polynomial of degree $\leq \dim M -2$,
so it must be a constant ($\dim M \leq 2$). 
Thus $\length (J^jM/I^jM)$ is a constant for $j \gg 0$.

For the second part, let us observe that
\begin{align*}
\length \Big(&\frac{J^jM}{I^jM}\Big) - \length \Big(\frac{\fm^{i+1}J^jM+J^{j+1}M}{\fm^{i+1}I^jM+I^{j+1}M}\Big)\\
& =\length \Big(\frac{J^jM}{\fm^{i+1}J^jM+J^{j+1}M}\Big) - \length \Big(\frac{I^jM}{\fm^{i+1}I^jM+I^{j+1}M}\Big)\\
 &=[j_1^1(J,M)-j_1^1(I,M)]i+ [j_0(J,M)-j_0(I,M)]j + j_1^2(J,M) - j_1^2(I,M).
\end{align*}
By \ref{am-prop}, it follows
that $$j_0(J,M)=j_0(I,M) \text{ and  } j_1^1(J,M)=j_1^1(I,M),$$
 and  therefore the last expression is a constant. 
Using the first part we can now conclude  the second part. 

By  Lemma \ref{fin_lemma} and Lemma \ref{dade_lemma}, we have 
\begin{align*}
\length (J^{j+1}M/I^{j+1}M)&=\length
\Big(G_\fm(J^{j+1}M)/G_\fm(I^{j+1}M) \Big)\\&=\length \Big( \bigoplus_{k
  \geq 0}
\frac{(J^{j+1}M \cap \fm^kM) + \fm ^{k+1}M}{(I^{j+1}M \cap \fm^kM) +
  \fm^{k+1}M} \Big)\\
&=\length \Big( \bigoplus_{k=0}^{t}
\frac{(J^{j+1}M \cap \fm^kM) + \fm ^{k+1}M}{(I^{j+1}M \cap \fm^kM) +
  \fm^{k+1}M} \Big)
\end{align*}
  for some fixed integer $t$ independent of $j$ (by part (1) we can do
  this). Similarly,
\begin{align*}
\length \Big(\frac{\fm^{i+1}J^jM + J^{j+1}M}{\fm^{i+1}I^jM+ I^{j+1}M}\Big)&= \length
\Big(\frac{G_\fm(\fm^{i+1}J^jM + J^{j+1}M)}{G_\fm(\fm^{i+1}I^jM+
  I^{j+1}M)} \Big)\\ 
&=\length \Big( \bigoplus_{k=0}^{s}
\frac{\big((\fm^{i+1}J^jM + J^{j+1}M) \cap \fm^kM \big) + \fm
  ^{k+1}M}{\big((\fm^{i+1}I^jM+ I^{j+1}M) \cap \fm^kM \big) +
  \fm^{k+1}M} \Big),
\end{align*}
 for some fixed integer $s$ independent of $i$ and $j$ (we use here
 the second part of the statement). We may assume $s=t$. On the other hand, for $ i \geq t$, 
$$ (\fm^{i+1}J^jM + J^{j+1}M) \cap \fm^kM = \fm^{i+1} J^jM + (\fm^kM
\cap J^{j+1}M)$$ and 
$$ (\fm^{i+1}I^jM + I^{j+1}M) \cap \fm^kM = \fm^{i+1} I^jM + (\fm^kM
\cap I^{j+1}M).$$
This implies that
$$ \big((\fm^{i+1}J^jM + J^{j+1}M) \cap \fm^kM \big)+ \fm^{k+1}M =(\fm^kM
\cap J^{j+1}M) +\fm^{k+1}M$$ 
and 
$$ \big((\fm^{i+1}I^jM + I^{j+1}M) \cap \fm^kM \big)+ \fm^{k+1}M =(\fm^kM
\cap I^{j+1}M) +\fm^{k+1}M.$$ 

 We then get
\begin{align*}
\length \Big(\frac{\fm^{i+1}J^jM + J^{j+1}M}{\fm^{i+1}I^jM+
  I^{j+1}M}\Big)&=\length \Big( \bigoplus_{k=0}^{t}
\frac{(J^{j+1}M \cap \fm^kM) + \fm ^{k+1}M}{(I^{j+1}M \cap \fm^kM) +
  \fm^{k+1}M} \Big)\\&=\length (J^{j+1}M/I^{j+1}M)\\
&=\length (J^{j}M/I^{j}M),
\end{align*}
where the last equality follows from part (1).
\end{proof}
\begin{Lemma} Let $(A,\fm)$ be a local ring, let $I \subseteq J$ be two ideals
  in $A$, and let $M$ be a finitely generated $A$-module.  Let $k$ be
  a positive integer.
\begin{itemize}
\item[\rm 1)] If $I$ is a reduction of $(J,M)$,  then $
  j_0(J,M)=j_0(J,I^kM)$.
\item[\rm 2)] If $I$ is a
  reduction of $(J,M)$, then $I$ is a reduction of $(J,I^kM)$. 
\item[\rm 3)] Assume that $\dim M/IM < \dim M$ and that $M$ is
  equidimensional. If $I$ is a reduction of $(J,I^kM)$, then $I$ is a
  reduction of $(J,M)$ .
\item[\rm 4)] If $I$ is a reduction of $(J,M)$, then $j_1(I,I^kM)=j_1(J,I^kM)$
  implies that $j_1(I,M)=j_1(J,M)$.
\item[\rm 5)] $J \subseteq I_{\{1\}}^M$ if and only if $J \subseteq
  I_{\{1\}}^{I^kM}$.
\end{itemize}
\end{Lemma}
\begin{proof}(1) $I$ is a reduction of $(J,M)$, so there exists a
  positive integer $n$ such that $IJ^nM=J^{n+1}M$. So
  for $j \gg 0$,
\begin{equation}\label{id}
\length\Big(\frac{J^jI^kM}{\fm^{i+1}J^jI^kM +
  J^{j+1}I^kM}\Big)=\length\Big(\frac{J^{j+k}M}{\fm^{i+1}J^{j+k}M +
  J^{j+1+k}M}\Big),
\end{equation}
which implies that $j_0(J,M)=j_0(J,I^kM)$.

(4) also follows from (\ref{id}).

(2) is obvious.

(3) Let $\overline A=A/\Ann M$, $\bar{I}=I\overline{A}$, and $\bar{J}=J\overline{A}$.
 Then $\overline A$ is an equidimensional ring
and $\bar{I}$ is not contained in any minimal prime ideal of $\overline A$. 
Since $I$ is a reduction of $(J,I^kM)$, there exists a positive
integer $n$ such that $IJ^nI^kM=J^{n+1}I^kM$. By
the determinant trick, it follows that
$\bar{I}\bar{J}^n\bar{I}^k=\bar{I}^{k+1}\bar{J}^n$ is a reduction of
$\bar{J}^{n+1}\bar{I}^k$, so there exists $l$ such that 
$$\bar{I}^{k+1}\bar{J}^n(\bar{J}^{n+1}\bar{I}^k)^l=(\bar{J}^{n+1}\bar{I}^k)^{l+1}.$$
Set $s=kl+k$, $t=nl+n+l$ so that the above equality can be written
$$\bar{I}(\bar{I}^s\bar{J}^t)=\bar{J}(\bar{I}^s\bar{J}^t).$$

We claim that this implies that  $\bar{I}$ is a reduction of
$\bar{J}$. It is enough to show this after we mod out an arbitrary minimal prime ideal of
$\overline{A}$, and since $\bar{I}$ is not contained in any minimal
prime ideal of $\overline{A}$, we may therefore assume that $\overline{A}$ is a
domain and $\bar{I}$, $\bar{J}$ are nonzero ideals. Using again the
determinant trick, we get  $\bar{I}$ is a reduction of $\bar{J}$
($\bar{I}^s\bar{J}^t \not=0$),
which implies that $I$ is a reduction of $(J,M)$.

(5) Denote $K=\bigoplus_{n \geq 0} JI^nM/I^{n+1}M$ and $ L=\bigoplus_{n \geq
   0} JI^{n+k}M/I^{n+1+k}M$. 
It is clear that $\dim K=\dim L$. On the
   other hand, $J \subseteq I_{\{1\}}^M$ if and only if $\dim K < \dim
   G_I(M)$, and $J \subseteq I_{\{1\}}^{I^kM}$ if and only if $\dim L < \dim
   G_I(M)$.
\end{proof}

The following proposition shows that the first two generalized Hilbert
coefficients are the same up to the first coefficient
ideal.
\begin{Proposition} \label{implication1}Let $(A,\fm)$ be a local ring,
  let $M$ be a formally
  equidimensional $A$-module, and let $I$ be an ideal of $A$ with $\dim M/IM <
  \dim M$. If $I \subseteq J
  \subseteq I_ {\{1\}}^M$, then
  $$j_i(I,M)=j_i(J,M) \quad \text{for } i=0,1.$$
\end{Proposition}

\begin{proof} We may assume that $A$ is complete and that $M$ is
  equidimensional.

If $\dim M=1$, then the conclusion follows from Shah's
  result (in its version for modules). Indeed, we can replace $A$ by   $A/\Ann
  M$, and then the ideals $I$ and $J$ are primary to the maximal ideal
  of $A/\Ann M$.

  If $\dim M=2$,  from Proposition \ref{lowdim} part (3) it follows
  that for $i,j \gg0$ we have following equality of polynomial
  functions of degree one:
$$\length(I^jM/\fm^{i+1}I^jM+I^{j+1}M)=\length(J^jM/\fm^{i+1}J^jM+J^{j+1}M).$$
By Remark \ref{simple}, it follows that  $j_i(I,M)=j_i(J,M) \quad \text{for } i=0,1.$

  Assume $\dim M \geq 3$. If $\depth_I M=0$, replacing $M$ by $I^kM$
  for $k$ big enough, we may assume $\depth_I M > 0$ (the previous
  proposition shows that the hypotheses are preserved). 

  By 
  Proposition \ref{sh_lem1}, there exists an integer $n \geq 1$ and an
  element $ a
  \in I^n \setminus I^{n+1}$, with $a^{*}$ part of a system of parameters  of $G_I(M)$, such that 
  $I_{\{1\}}^M=(I^{n+1 }:aM)$. Since $I$ is a reduction of $J$
  (see Proposition \ref{red}),  we
  can choose $x \in I \setminus \fm J$ superficial element for
  $(J,M)$ ($I'$ and $J'$ have the same radical in $G_\fm(G_J(A))$
  ). By taking a sufficiently general element,  we may  also assume  that
  $x$  is a superficial
  element for $(I,M)$,  $a^{*},x^{*}$ are part of a system of
  parameters  of
  $G_\fm(G_I(A))$, and  $x$  is   a nonzerodivisor on
  $M$ ($\depth_I M > 0$).

  Denote $\overline{M}=M/xM$. By the choice of $x$ it follows that $I
  \subseteq J \subseteq I_ {\{1\}}^{\overline{M}}$. Indeed, if $y \in
  J$, then $ya \overline{M} \subseteq I^{n+1}\overline{M}$. But $x^{*}$ and $a^{*}$
  are part of a system of
  parameters of $G_I(M)$, so $a^{*}$ is part of a system of
  parameters
  of $G_I(\overline{M}) \cong G_I(M)/x^{*}G_I(M) $. Then  $\bar y \in
  I_ {\{1\}} ^{\overline{M}}$ and the induction hypothesis gives
  $j_i(I,\overline{M})=j_i(J,\overline{M})$ for $i=0,1.$
  Using Proposition \ref{dif} we now obtain $j_i(I,M)=j_i(J,M)$
  for  $i=0,1$.

 Note that we
 cannot prove the $2$-dimensional case by reducing the problem to the
 $1$-dimensional case.
 The polynomial that gives $\length(I,M)
 (I^jM/\fm^iI^jM+I^{j+1}M)$ for $j \gg 0$ has
 the form $j_1^1(I,M)i+j_0(I,M)j+j_1^2(I,M)$. By reducing the dimension one
 more time we would
 loose the coefficients $j_1^1(I,M)$ and $j_1^2(I,M)$.    
 \end{proof}

We can now prove the theorem stated in the introduction.
\begin{Theorem}\label{mainth} Let $\ringA$ be a local ring, let $M$ be
a formally
  equidimensional $A$-module, and let $I \subseteq J$ be two ideals of
  $A$ with $\dim M/IM < \dim M$. The following are equivalent:
\begin{itemize}
\item[\rm 1)] $J \subseteq I_{\{1\}}^M$.
\item[\rm 2)] $j_i(I_\fp,M_\fp)=j_i(J_\fp,M_\fp)$ for $i=0,1$ and every
  $\fp \in \Spec(A)$.
\end{itemize}
\end{Theorem}

\begin{proof} The proof of the  case $\dim M =2$ is the crucial
part of the argument.  Then we can use  an induction argument similar 
to the one
used by Flenner and Manaresi in the proof of their theorem (see the
introduction).    

If $\dim M=1$, using the same argument used in the proof of the
previous theorem, we can reduce the problem to the $\fm$-primary case and Shah's result proves both implications.

As usual, we may assume that $\ringA$ is a complete local
  ring and $M$ is equidimensional. We will prove that for every prime
  ideal $\fp$, $J_\fp \subseteq {(I_\fp)}_{\{1\}}^{M_\fp}$ and
  the implication $(1) \Longrightarrow (2)$ will follow from
  Proposition \ref{implication1}. 

Let $N=\bigoplus_{n \geq 1} JI^n/I^{n+1}$. Since
  $J \subseteq
  I_{\{1\}}^M$, by Remark \ref{l_dim}, we have  $\dim N < \dim G_I(M)=\dim
  M$. Let $N'= \bigoplus_{n \geq 1}
  J_\fp(I_\fp)^n/(I_\fp)^{n+1}= U^{-1}N$, where $U=G_I^0(A) \setminus (\fp/I)$
  is a multiplicatively closed subset of $G_I(A)$. Since $G_I(M)$ is
  equidimensional, we get
  $\dim N' < \dim U^{-1}G_I(M)=\dim G_{I_\fp}(M_\fp)$, i.e. $J_\fp
  \subseteq (I_\fp)_{\{1\}}^{M_\fp}$.

We prove the converse by induction on $d=\dim M$. First assume $\dim M=2$. We can also assume
that $M$ is faithful, so $\dim A=2$. 
Since $j_i(I,M)=j_i(J,M)$ for $i=0,1$ there exist $i_0,j_0$ such that
for $i \geq i_0$ and $j \geq j_0$ 
\begin{equation}\label{in}
\length(I^jM/\fm^{i+1}I^jM+I^{j+1}M)=\length(J^jM/\fm^{i+1}J^jM+J^{j+1}M).
\end{equation}
Let $\fp \in
\Spec(A) \setminus \{\fm\}$, so by hypothesis $j_i(I_\fp,M_\fp)=j_i(J_\fp,M_\fp)$ for
$i=0,1$. But $\dim A_\fp=1$, so  $I_\fp$ and $J_\fp$ are primary to
the maximal ideal. Applying the theory of first coefficient ideals
for $\fm$-primary ideals (in a version for modules) we get
$\length(J^j_\fp M_\fp/I^j_\fp M_\fp)=0$ for $j \gg 0$ 
(it is bounded above by
a polynomial of degree $\dim A_\fp -2=-1$). There are only finitely
many elements in $\Spec(A) \setminus \{\fm\}$ that contain $I$, so there exists $r \geq j_0$ such
that for all $\fp \in \Spec(A) \setminus \{\fm\}$ and $j \geq r$ 
we have $\length(J^j_\fp M_\fp/I^j_\fp M_\fp)=0$, and this   implies that
$\length (J^jM/I^jM) < \infty$ for $j \geq r$. Choose $c \geq i_0$
such that 
$\fm^i J^rM \subseteq I^rM$ for $i \geq c$.

We are now using an argument similar to the one given  in the proof of
Proposition \ref{lowdim}.

For $i \geq c$, we have 
\begin{align*}
(*)&\length \Big(\frac{J^rM}{I^rM}\Big) - \length \Big(\frac{\fm^{i+1}J^rM+J^{r+1}M}{\fm^{i+1}I^rM+I^{r+1}M}\Big)\\
& =\length \Big(\frac{J^rM}{\fm^{i+1}J^rM+J^{r+1}M}\Big) - \length \Big(\frac{I^rM}{\fm^{i+1}I^rM+I^{r+1}M}\Big)\\
 &=[j_1^1(J,M)-j_1^1(I,M)]i+ [j_0(J,M)-j_0(I,M)]j + j_1^2(J,M) -
 j_1^2(I,M)\\
&=0.
\end{align*}
where the last equality follows from hypothesis.

Then, by Lemma \ref{fin_lemma} and Lemma \ref{dade_lemma}, we have 
\begin{align*}
\length (J^{r+1}M/I^{r+1}M)&=\length
\Big(G_\fm(J^{r+1}M)/G_\fm(I^{r+1}M) \Big)\\&=\length \Big( \bigoplus_{k
  \geq 0}
\frac{(J^{r+1}M \cap \fm^kM) + \fm ^{k+1}M}{(I^{r+1}M \cap \fm^kM) +
  \fm^{k+1}M} \Big)\\
&=\length \Big( \bigoplus_{k=0}^{t}
\frac{(J^{r+1}M \cap \fm^kM) + \fm ^{k+1}M}{(I^{r+1}M \cap \fm^kM) +
  \fm^{k+1}M} \Big)
\end{align*}
for some fixed integer $t$ ($r$ is fixed). 

Using $(*)$ we obtain that  for $i \geq c$ 
\begin{align*} 
(**)\length(J^rM/I^rM)&=\length \Big(\frac{\fm^{i+1}J^rM + J^{r+1}M}{\fm^{i+1}I^jM+ I^{r+1}M}\Big)\\&= \length
\Big(\frac{G_\fm(\fm^{i+1}J^rM + J^{r+1}M)}{G_\fm(\fm^{i+1}I^rM+
  I^{r+1}M)} \Big)\\ 
&=\length \Big( \bigoplus_{k=0}^{s}
\frac{\big((\fm^{i+1}J^rM + J^{j+1}M) \cap \fm^kM \big)+ \fm ^{k+1}M}{\big((\fm^{i+1}I^rM+ I^{r+1}M) \cap \fm^kM \big)+
  \fm^{k+1}M} \Big)
\end{align*}
 for some fixed integer $s$ independent of $i$ ($r$ is fixed). 
We may assume $s=t$.
But for  $ i \geq c+t$ and $0 \leq k \leq t$ we  have 
$$ (\fm^{i+1}J^jM + J^{r+1}M) \cap \fm^kM = \fm^{i+1} J^rM + (\fm^kM
\cap J^{r+1}M)$$ 
and the similar equality with $I$ instead of $J$. 
This implies that for $ i \geq c+t$
$$ \big((\fm^{i+1}J^jM + J^{r+1}M) \cap \fm^kM \big)+ \fm^{k+1}M =(\fm^kM
\cap J^{r+1}M) +\fm^{k+1}M$$ 
and the similar equality for $I$.
 
Using the above observations and $(**)$ we have that for $i \geq c+t$ 
\begin{align*}
\length(J^rM/I^rM)&= \length \Big(\frac{\fm^{i+1}J^rM + J^{r+1}M}{\fm^{i+1}I^rM+
  I^{r+1}M}\Big)\\
&=\length \Big( \bigoplus_{k=0}^{t}
\frac{(J^{r+1}M \cap \fm^kM) + \fm ^{k+1}M}{(I^{r+1}M \cap \fm^kM) +
  \fm^{k+1}M} \Big)\\
&=\length (J^{r+1}M/I^{r+1}M).
\end{align*}

Repeating the argument we conclude that
$\length(J^jM/I^jM)$ is constant for $j \geq r$.

But this implies that there exists $l$ such that $\fm^l J^jM \subseteq
I^jM$ for all $j \geq r$. Indeed, if $L$ is module of finite length,
say $l$, over a local ring, then $\fm^lL=0$. 

So $\fm^iJ^jM \subseteq I^jM$ for all $i \geq l$ and $j \geq r$, hence  
$ \fm^iJI^jM \subseteq \fm^iJ^{j+1}M \subseteq I^{j+1}M$ for all $i
\geq l$ and $j \geq r$.
Then $$\length \Big(\frac{\fm^iJI^jM + I^{j+1}M}{\fm^{i+1}JI^jM + I^{j+1}M}\Big)=0
  \quad \text{for } i \geq l, j \geq r ,$$
which implies that $\dim_{G_\fm(G)}G_\fm(N) \leq 1$, where $N=\bigoplus_{n \geq
  0}JI^nM/I^{n+1}$. This means that $\dim N \leq 1=\dim G_I(M)-1$, and
  by Remark \ref{l_dim} we get $J \subseteq I_{\{1\}}^M$.

We now assume that $\dim M \geq 3$.
Replacing $M$ by $I^kM$ for a suitable $k$, we may assume
  that $\depth _I(M) > 0$. 

By the theorem of Flenner and Manaresi (\cite[Theorem 3.3]{FM}; see
also the introduction), 
$I$ is a reduction of $(J,M)$. Then, for  a
sufficiently general element
$x$ in $I$, we have $j_i(I_\fp,M_\fp)=j_i(I_\fp,
\overline{M}_\fp)$ and $j_i(J_\fp,M_\fp)=j_i(J_\fp,
\overline{M}_\fp)$ for $i=0,1$, where  $\overline{M}=M/xM$. By the induction
hypothesis, we get $J \subseteq I_{\{1\}}^{\overline{M}}$. We still
have to prove that   $J \subseteq I_{\{1\}}^{M}$. 

Let $K=\bigoplus_{n \geq 0} JI^nM/I^{n+1}M$ and 
$L=\bigoplus_{n \geq 0} (JI^nM+xM)/(I^{n+1}M+xM)$. Since $J \subseteq
I_{\{1\}}^{\overline{M}}$, we have $\dim L < \dim G_I(M/xM)=\dim M/xM=\dim M -1$
(we can choose $x$ to be a nonzerodivisor on $M$). 

For technical reasons (see Proposition \ref{dif}), we will prefer this interpretation of the Generalized Hilbert coefficients.

Consider the exact sequence 
\begin{equation}\label{exact}
0\to U{\to} K\stackrel{\pi}{\to} L\to 0
\end{equation} 
where $U$ is the kernel of the canonical epimorphism $K\stackrel{\pi}{\to}L$.

We have
\begin{align*}
U_n&=\frac{JI^nM \cap(I^{n+1}M+xM)}{I^{n+1}M}= \frac{I^{n+1}+(JI^nM
  \cap xM)}{I^{n+1}M} \\
&\cong \frac{JI^nM \cap xM}{I^{n+1}M \cap
    xM}=\frac{x(JI^nM:x)}{x(I^{n+1}M:x)} \cong
  \frac{JI^nM:x}{I^{n+1}M:x}.
\end{align*}
On the other hand, for  $x \in I$ sufficiently general and $n \gg 0$,  $(JI^{n+1}M:x)=JI^nM$ and
$(I^{n+1}M:x)=I^nM$.

So, for $n$ large enough, $K_{n-1} \cong U_{n}$ (isomorphism induced by
the multiplication by $x$), and then the exact sequence
\ref{exact} implies that
$\dim L=\dim K - 1$.  
Since $\dim L < \dim M -1$
we have $\dim K < \dim M$, i.e. $J \subseteq I_{\{1\}}^M$.
\end{proof}

We now sketch a proof of Proposition \ref{lat}. 

\begin{proof}[Proof of Proposition \ref{lat}\/]
Note that by \ref{am-prop}
we have $c_d(I,M)=c_d(J,M)$, so all we need to prove is that $c_i(I,M)=c_i(J,M)$ for $i=0,\ldots,d-1$. We use induction on $d=\dim M$.

If $d=0,1$, the conclusion follows
immediately from  \ref{am-prop}. Replacing $M$ by $I^kM$ for a suitable
$k$, we may assume that $\depth_I(M) > 0$. Indeed, $I$ is a reduction of $(J,M)$, so there exists a
  positive integer $n$ such that $IJ^nM=J^{n+1}M$. Then
  for $j \gg 0$,
\begin{equation}
\length\Big(\frac{J^jI^kM}{\fm^{i+1}J^jI^kM +
  J^{j+1}I^kM}\Big)=\length\Big(\frac{J^{j+k}M}{\fm^{i+1}J^{j+k}M +
  J^{j+1+k}M}\Big),
\end{equation}
which implies that $c_i(J,M)=c_i(J,I^kM)$ for $i=0,1,\ldots,d-1$. 

Choose $x \in I$  a nonzerodivisor on $M$ which is a superficial
element for $(I,M)$ and $(J,M)$. By Proposition \ref{dif}, we have
$c_i(I,M)=c_i(\overline{I},\overline{M})$ and $c_i(J,M)=c_i(\overline{J},\overline{M})$
for $i=0,\ldots,d-1$, where for an $A$ module $L$ we denote
$\overline{L}=L/xL$. 
The induction  hypothesis implies that $c_i(I,M)=c_i(J,M)$ for $i=0,\ldots,d-1$.
\end{proof}

\begin{Example} Let $A=k[x,y,z]$ be the ring of polynomials in three variables over the field $k$,  and let $\fm=(x,y,z)$ be the maximal homogeneous ideal. As in the local case, one can define the generalized Hilbert coefficients and the first coefficient ideal associated to an ideal. 

 Let $I=(x^5,y^3,xyz^2)$ and  let $J=(x^5,y^3,xyz^2,x^4y^2)$. Note that both ideals have  height $2$ and analytic spread $3$.   A computation with Macaulay 2 \cite{Mac} shows that $j_0(I)=j_0(J)=30$, $j_1(I)=j_1(J)=(8,-32)$, $j_2(I)=(0,-1,5)$, $j_2(J)=(0,-1,3)$. 

In fact, using the method described in \cite[Proposition 3.2]{C}, one can show  that $J=I_{\{1\}}$, hence the equality of the first two generalized Hilbert coefficients.

\end{Example}
\begin{ack} The author thanks   Craig Huneke for
  valuable discussions concerning the material
  of this paper.
\end{ack}

\end{document}